\documentclass[12pt]{amsart}
\usepackage{amssymb}
\input amssym.def
\usepackage{amsmath, amsfonts,hyperref}
\usepackage{amscd}
\usepackage[mathscr]{eucal}

\usepackage{xcolor}
\hypersetup{
    colorlinks,
    linkcolor={red!50!black},
    citecolor={blue!50!black},
    urlcolor={blue!80!black}
}

\newcommand{\Z}{{\mathbb Z}}
\newcommand{\Q}{{\mathbb Q}}
\newcommand{\C}{{\mathbb C}}

\newcommand{\cT}{{\mathcal T}}
\newcommand{\cZ}{{\mathcal Z}}
\newcommand{\z}{\zeta}

\newtheorem*{conj}{Conjecture}

\DeclareFontFamily{U}{wncy}{}
\DeclareFontShape{U}{wncy}{m}{n}{<->wncyr10}{}
\DeclareSymbolFont{mcy}{U}{wncy}{m}{n}
\DeclareMathSymbol{\Sh}{\mathord}{mcy}{"58}

\begin{document}

\title{A conjecture about multiple $t$-values}

\author{Biswajyoti Saha}

\address{Biswajyoti Saha\\ \newline
School of Mathematics, Tata Institute of Fundamental Research, Dr. Homi Bhabha Road,
Navy Nagar, Mumbai 400005, India.}
\email{biswa@math.tifr.res.in}



\begin{abstract}
For positive integers $a_1,\ldots,a_r$ with $a_1 \ge 2$, the multiple
$t$-value $t(a_1,\ldots,a_r)$ is defined by the series
$\sum\limits_{n_1 > \ldots > n_r > 0 \atop n_i \text{ odd}} n_1^{-a_1} \cdots n_r^{-a_r}$.
For an integer $k \ge 2$, the dimension of the $\Q$-vector space generated by all the
multiple $t$-values of weight $k$ has been predicted by Hoffman \cite{MEH} to be
the $k$-th Fibonacci number. In this short note we give a conjectural basis
of this vector space.
\end{abstract}

\date{\today}

\maketitle

Let $r \ge 1$ be an integer and
$$
U_r:=\{(s_1,\ldots,s_r) \in \C^r : \Re(s_1 + \cdots +s_i)>i \ \text{for} \ 1 \le i \le r \}.
$$
Following \cite{MEH}, we define the multiple $t$-function of depth $r$ by
the infinite series
$$
t(s_1,\ldots,s_r):=\sum_{n_1 > \ldots > n_r > 0 \atop n_i \text{ odd}} n_1^{-s_1} \cdots n_r^{-s_r}
$$
for $(s_1,\ldots,s_r) \in U_r$. We further have the equality
$$
t(s_1,\ldots,s_r)=2^{-(s_1+\cdots+s_r)} \z(s_1,\ldots,s_r; -\frac{1}{2}, \ldots, -\frac{1}{2})
$$
where $\z(s_1,\ldots,s_r; -\frac{1}{2}, \ldots, -\frac{1}{2})$ is the multiple Hurwitz zeta function
$$
\sum_{n_1 > \ldots > n_r > 0} (n_1-1/2)^{-s_1} \cdots (n_r-1/2)^{-s_r}
$$
defined on $U_r$. The meromorphic continuation of the multiple Hurwitz zeta functions is well known
(see \cite{AI}). The exact set of polar hyperplanes of the multiple Hurwitz zeta functions is given
in \cite{GS}. Thus we get that the multiple $t$-function $t(s_1,\ldots,s_r)$ has a meromorphic
extension to $\C^r$ with polar singularities along the hyperplanes
$$
H_{1,0}, H_{2,1}, H_{2,2k} \text{ and } H_{i,k} 
\text{ for all } k\ge 0, 3 \le i \le r,
$$
where for $i \ge 1$ and $k \ge 0$, $H_{i,k}$ denotes the hyperplane  defined by
$s_1+\cdots+s_i=i-k$.

The special values $t(a_1,\ldots,a_r)$ for integers $a_i \ge 1$ and $a_1 \ge 2$
are the central object of study in \cite{MEH}. Hoffman \cite{MEH} studied these numbers
in the spirit of the multiple zeta values
$\z(a_1,\ldots,a_r):=\sum_{n_1 > \ldots > n_r > 0} n_1^{-a_1} \cdots n_r^{-a_r}$.
Hoffman noted various striking similarities of these numbers with the multiple
zeta values. He also pointed out some major differences. The sum $a_1+\cdots+a_r$
is called the weight of the values $\z(a_1,\ldots,a_r)$ and $t(a_1,\ldots,a_r)$.

Let $k \ge 2$ be an integer. Zagier considered the $\Q$-vector space generated by all the
multiple zeta values of weight $k$, which we denote by $\cZ_k$, i.e.
$$
\cZ_k := \Q \langle \z(a_1,\ldots,a_r) : a_1+\cdots+a_r=k \rangle.
$$
Zagier \cite{DZ} predicted that $\dim_\Q \cZ_k$ satisfies the recurrence relation
$$
d_k=d_{k-2}+d_{k-3}
$$
for $k \ge 5$ with the initial values $d_2=d_3=d_4=1$. Deligne-Goncharov \cite{DG} and
Terasoma \cite{TT} independently proved that $\dim_\Q \cZ_k \le d_k$.

In \cite{MEH}, Hoffman considered the analogous $\Q$ vector space generated by all the
multiple $t$-values of weight $k$, which we denote by $\cT_k$, i.e.
$$
\cT_k := \Q \langle t(a_1,\ldots,a_r) : a_1+\cdots+a_r=k \rangle.
$$
Finding various relations among multiple $t$-values of weight up to $7$, he
conjectured that $\dim_\Q \cT_k$ satisfies the recurrence relation
$$
f_k=f_{k-1}+f_{k-2}
$$
for $k \ge 4$ with the initial values $f_2=1,f_3=2$. Note that this is the recurrence
relation of the Fibonacci numbers $F_k$.

Further, for $\cZ_k$, Hoffman \cite{MEH1} conjectured that the set
$$
B_k:=\{\z(a_1,\ldots,a_r) : a_1+\cdots+a_r=k, a_i \in \{2,3\}  \text{ and } 1\le r \le k-1\}
$$
forms a $\Q$-basis. It is consistent with the above conjecture of Zagier i.e.
$|B_k| \le d_k$. Brown \cite{FB} proved that $B_k$ is in fact a generating set.

In this short note, we give a conjectural basis for $\cT_k$ for $k \ge 2$. Let 
$$
C_k:=\{t(a_1+1,a_2,\ldots,a_r) : a_1+\cdots+a_r=k-1, a_i \in \{1,2\} \text{ and } 1\le r \le k-1\},
$$
for any $k \ge 2$. For example,
\begin{align*}
& C_3=\{t(2,1),t(3)\}, \ C_4=\{t(2,1,1),t(2,2),t(3,1)\} \text{ and}\\
& C_5=\{t(2,1,1,1), t(2,1,2), t(2,2,1), t(3,1,1), t(3,2)\}.
\end{align*}
It is not difficult to see that the cardinality of the set 
$$
\{(a_1,\ldots,a_r) : a_1+\cdots+a_r=k-1, a_i \in \{1,2\} \text{ and } 1\le r \le k-1 \}
$$
is the $k$-th Fibonacci number $F_k$. We now conjecture the following.

\begin{conj}
The following set of weight $k$ multiple $t$-values
$$
\{t(a_1+1,a_2,\ldots,a_r) : a_1+\cdots+a_r=k-1,  a_i \in \{1,2\} \text{ and } 1\le r \le k-1 \}
$$
forms a $\Q$-basis of $\cT_k$ for $k \ge 2$.
\end{conj}

This conjecture is in accordance with various relations among multiple
$t$-values of weight up to $7$ found by Hoffman \cite[Appendix]{MEH}, i.e. for any
weight $2 \le k \le 7$, it follows from \cite[Appendix]{MEH} that $C_k$ generates
$\cT_k$. As an example, we give explicit calculations for weight $5$.

Note that $C_5=\{t(2,1,1,1), t(2,1,2), t(2,2,1), t(3,1,1), t(3,2)\}$.
We need the following four explicit relations from \cite[Appendix]{MEH}:
\begin{align}\label{5-2,3}
t(3,2)= -\frac{1}{2} t(5) + \frac{3}{7} t(2)t(3), \ \
t(2,1,2)= \frac{3}{4} t(5) - \frac{1}{2} t(2)t(3),
\end{align}
and
\begin{align}\label{4,1}
t(4,1)= -\frac{1}{2} t(5) - \frac{1}{7} t(2)t(3) + t(4) \log 2, \ \
t(2,2,1)= \frac{1}{8} t(5) - \frac{3}{14} t(2)t(3) + \frac{1}{4}t(4) \log 2.
\end{align}
From \eqref{5-2,3} (and the fact that $t(2)t(3)=t(2,3)+t(3,2)+t(5)$), one can
write $t(3,2)$ and $t(2,1,2)$ as $\Q$-linear combination of $t(5)$ and $t(2,3)$ as follows:
$$
t(3,2)= -\frac{1}{8} t(5) + \frac{3}{4} t(2,3), \ \
t(2,1,2)= \frac{5}{16} t(5) - \frac{7}{8} t(2,3).
$$
This system of equations is invertible. Therefore $t(5)$ and $t(2,3)$
can be written as $\Q$-linear combination of $t(3,2)$ and $t(2,1,2)$.
Now from \eqref{4,1}, one can write $t(4,1)$ as a $\Q$-linear combination
of $t(3,2), t(2,2,1)$ and $t(2,1,2)$.

\end{document}